\documentclass[12pt]{article}
\usepackage{amsfonts,amssymb,amsthm}
\begin{document}
\unitlength=1mm
\newtheorem{remark}{Remark}
\newtheorem{theorem}{Theorem}
\newtheorem{proposition}{Proposition}
\newtheorem{plain}{Definition}
\newtheorem{cor}{Corollary}
\newcommand{\Mt}{\mathop{\rm Mult}\nolimits}
\title{Quantum seaweed algebras and quantization of\\
 affine Cremmer--Gervais $r$-matrices}
\author{M.E. Samsonov$^{1)}$, A.A. Stolin$^{2)}$ and V.N. Tolstoy$^{3)}$
\\[10pt]
$^{1)}$\small{\sl St. Petersburg State University, Institute of Physics,}
\\[-2pt]
\small{\sl198504 St. Petersburg, Russia; e-mail: samsonov@pink.phys.spbu.ru}
\\[5pt]
$^{2)}$\small{\sl Department of Mathematics, University of G\"oteborg,}
\\[-2pt]
\small{\sl SE-41296 G\"oteborg, Sweden; e-mail: astolin@math.chalmers.se}
\\[5pt]
$^{3)}$\small{\sl Institute of Nuclear Physics, Moscow State University},
\\[-2pt]
\small{\sl119992 Moscow, Russia;  e-mail: tolstoy@nucl-th.sinp.msu.ru}}
\date{}

\maketitle \abstract{We propose a method of quantization of certain Lie
bialgebra structures on the polynomial Lie algebras related to
quasi-trigonometric solutions of the classical Yang--Baxter equation. The
method is based on an affine realization of certain seaweed algebras and their
quantum analogues. We also propose a method of $\omega$-affinization, which
enables us to quantize rational $r$-matrices of $\mathfrak{sl}(3)$.}

\section{Introduction}
The aim of this paper is to propose a method of quantizating certain
Lie bialgebras structures on polynomial Lie algebras. In the beginning of
90-th Drinfeld in \cite{D} posed the following problem: can any Lie bialgebra
be quantized? The problem was solved by Etingof and Kazhdan and the answer
was positive. However, another problem of finding explicit quantization
formulas remains open.

First results in this direction were obtained in \cite{KM}, \cite{ESS}  and
\cite{IO}, where the authors quantized the so-called Belavin--Drinfeld list,
the list of all quasi-triangular Lie bialgebra structures on finite dimensional
simple Lie algebras. It should be noticed that first infinite-dimensional
cases were considered in \cite{KST}. However, a real break-through in the
infinite dimensional case came in \cite{KST1} and \cite{Sam1}, where deformed
versions of Yangians $Y(\mathfrak{sl}_{N}^{})$ and quantum affine algebras
$U_{q} (\hat{\mathfrak{sl}}_{N}^{})$ were constructed for $N=2,3$. In the
present paper we solve the problem in the case $U_q (\hat{\mathfrak{sl}}_{N}^{})$.
Our solution of the problem is based on a $q$-version of the so-called
seaweed algebras. A seaweed subalgebra of $\mathfrak{sl}_{N}^{}$ is an
intersection of two parabolic subalgebras one of which containing the Borel
subalgebra $B^+$ and another $B^-$. The case when both parabolic subalgebras
are maximal was studied in \cite{S} in connection with the study of rational
solutions of the classical Yang-Baxter equation. In particular, a complete
answer to the question when such an algebra is Frobenius was obtained there.
Later in \cite{DK} it was found out when an arbitrary seaweed algebra is
Frobenius.

In the present paper we quantize a Lie bialgebra, which as a Lie algebra is
$\mathfrak{gl}_{N}^{}[u]$ ($\mathfrak{sl}_{N}^{}[u]$). Its coalgebra structure is defined by a
quasi-trigonometric solution of the classical Yang-Baxter equation.
Quasi-trigonometric solutions were introduced in \cite{KPST}. We remind briefly
some results obtained there.

For convenience we consider the case $\mathfrak{g}=\mathfrak{gl}_{N}^{}$
although the results are also valid for an arbitrary simple complex Lie algebra
$\mathfrak{g}$. Let $e_{ij}^{}$, $i,j=1,\ldots,N$, be the standard Cartan-Weyl
basis of $\mathfrak{gl}_{N}^{}$: $[e_{ij}^{},e_{kl}^{}]=\delta_{jk}^{}e_{il}^{}-
\delta_{il}^{}e_{kj}^{}$. The element $\mathcal{C}^{}_2:=\sum\limits_{i,j=1}^N
e_{ij}^{}e_{ji}^{}\in U(\mathfrak{gl}_{N}^{})$ is a $\mathfrak{gl}_{N}^{}$-scalar,
i.e. $[\mathcal{C}_{2}^{},x]=0$ for any $x\in \mathfrak{gl}_{N}^{}$,
and it is called the second order Casimir element. The element $\Omega:=
\frac{1}{2}\Big(\Delta(\mathcal{C}_{2}^{})-\mathcal{C}_{2}^{}\otimes1-
1\otimes\mathcal{C}_{2}^{}\Big)=\sum\limits_{i,j=1}^N e_{ij}^{}\otimes
e_{ji}^{}\subset U(\mathfrak{gl}_{N}^{})\otimes U(\mathfrak{gl}_{N}^{})$,
where $\Delta$ is a trivial co-product $\Delta(x)=x\otimes1+1\otimes x$
($\forall x\in \mathfrak{gl}_{N}^{}$), is called the Casimir two-tensor. The
two-tensor can be represented in the form $\Omega=\Omega_{+}+\Omega_{-}$,
where $\Omega_{+}=\frac{1}{2}\sum\limits_{1\leq i\leq N}e_{ii}^{}\otimes
e_{ii}^{}+\sum\limits_{1\leq i<j\leq N}e_{ij}^{}\otimes e_{ji}^{}$ and
$\Omega_{-}=\frac{1}{2}\sum\limits_{1\leq i\leq N}e_{ii}^{}\otimes e_{ii}^{}+
\sum\limits_{1\leq i<j\leq N}e_{ji}^{}\otimes e_{ij}^{}$. Note
that $(\omega\otimes\omega)(\Omega_{\pm})=\Omega_{\mp}$, where
$\omega$ is the Cartan automorphism: $\omega(e_{ij}^{})=-e_{ji}^{}$.

We say that a solution $X(u,v)$ of the classical Yang--Baxter equation is
{\it quasi-trigono\-metric} if it is of the form:
\begin{equation}\label{i1}
X(u,v)=\frac{u\,\Omega_{+}+v\,\Omega_{-}}{u-v}+p(u,v)~,
\end{equation}
where $p(u,v)$ is a non-zero polynomial with coefficients in
$\mathfrak{gl}_{N}^{}\otimes\mathfrak{gl}_{N}^{}$. If  $p(u,v)=0$ then $X(u,v)$
is the simplest (standard) trigonometric $r$-matrix.
Any quasi-trigonometric solution of the classical Yang-Baxter equation defines
a Lie bialgebra structure on $\mathfrak{gl}_{N}^{}[u]$ and the corresponding Lie
cobracket on $\mathfrak{gl}_{N}^{}[u]$ is given by the formula
\begin{equation}\label{i2}
\{A(u)\in\mathfrak{gl}_{N}^{}[u]\}\to \{[X(u,v),A(u)\otimes 1+1\otimes A(v)]\in
\mathfrak{gl}_{N}^{}[u]\otimes\mathfrak{gl}_{N}^{}[v]\}~.
\end{equation}
It was proved in \cite{KPST} that for
$\mathfrak{g}=\mathfrak{sl}_{N}^{}\,(\mathfrak{gl}_{N}^{})$ there is a one-to-one
correspondence between quasi-trigonometric $r$-matrices and Lagrangian subalgebras
of $\mathfrak{g}\oplus\mathfrak{g}$ transversal to a certain Lagrangian
subalgebra of $\mathfrak{g}\oplus\mathfrak{g}$ defined by a maximal parabolic
subalgebra of $\mathfrak{g}$. Here we mean the Lagrangian space with respect to the
following symmetric non-degenerate invariant bilinear form on
$\mathfrak{g}\oplus\mathfrak{g}$:
\begin{equation}\label{i3}
Q((a,b),(c,d))=K(a,c)-K(b,d)~,
\end{equation}
where $K$ is the Killing form and $a,b,c,d\in\mathfrak{g}$.

In their famous paper on the classical Yang--Baxter equation, Belavin and Drinfeld
listed all the trigonometric $r$-matrices. In case  $\mathfrak{sl}_{3}$ there exist 4
trigonometric $r$-matrices. Two of them relate to the quasi-triangular constant
$r$-matrices and can be quantized using methods from \cite{ESS}, \cite{IO} and
\cite{KM}. In our paper we explain how to to quantize one of the two remaining
trigonometric $r$-matrices found by Belavin and Drinfeld. We also notice that in
cases $\mathfrak{sl}_{2}$ and $\mathfrak{sl}_{3}$ our methods lead to quantization
of some rational $r$-matrices found in \cite{S}. As it was explained in \cite{Sam2},
quantization of the rational $r$-matrix for $\mathfrak{sl}_{2}$ has close relations
with the Rankin--Cohen brackets for modular forms (see \cite{CM}).

Our method is based on finding quantum twists, which are
various solutions of the so-called cocycle equation for a number of Hopf algebras:
\begin{equation}\label{i4}
F^{12}(\Delta\otimes{\rm id})(F)=F^{23}({\rm id}\otimes\Delta)(F)~.
\end{equation}
We are thankful to S.M. Khoroshkin and I. Pop for fruitful discussions and
suggestions.

\section{A quantum seaweed algebra and its affine realization}

It turns out that it is more convenient to use instead of the simple Lie
algebra $\mathfrak{sl}_{N}^{}$ its central extension $\mathfrak{gl}_{N}^{}$.
The polynomial affine Lie algebra $\mathfrak{gl}_{N}[u]$ is generated by
Cartan--Weyl basis $e_{ij}^{(n)}:=e_{ij}u^n$ ($i,j=1,2,\ldots N$,
$n=0,1,2,\ldots$) with the defining relations
\begin{equation}\label{pr1}
[e_{ij}^{(n)},\,e_{kl}^{(m)}]=\delta_{jk}e_{il}^{(n+m)}-\delta_{il}e_{kj}^{(n+m)}~.
\end{equation}
The total root system $\Sigma$ of the Lie algebra $\mathfrak{gl}_{N}[u]$ with
respect to an extended Cartan subalgebra generated by the Cartan elements
$e_{ii}$ ($i=1,2,\ldots, N$) and $d=u(\partial/\partial u)$ is given by
\begin{equation}\label{pr2}
\mathop{\Sigma}(\mathfrak{gl}_{N}[u])\,=\,\{\epsilon_i-\epsilon_j,\,n\delta+
\epsilon_i-\epsilon_j,\,n\delta~|~i\neq j;~i,j=1,2,\ldots,N;~n=1,2,\ldots\}~,
\end{equation}
where $\epsilon_i^{}$ ($i=1,2,\ldots,N$) is an orthonormal basis of a
$N$-dimensional Euclidean space $\mathbb R^{N}$: $(\epsilon_i^{},\epsilon_j^{})=
\delta_{ij}^{}$,  and we have following correspondence between 
root vectors and generators of $\mathfrak{gl}_{N}[u]$:
$e_{ij}^{(n)}=e_{n\delta+\epsilon_i-\epsilon_j}$ for $i\neq j$,
$n=0,1,2,\ldots$. It should be noted that the multiplicity of the
roots $n\delta$ is equal to $N$, $\Mt(n\delta)=N$, and the root vectors
$e_{ii}^{(n)}$ ($i=1,2,\ldots,N$) "split" this multiplicity.
We choose the following system of positive simple roots:
\begin{equation}\label{pr3}
\mathop{\Pi}(\mathfrak{gl}_{N}[u])\,=\,\{\alpha_i:=\epsilon_i-\epsilon_{i+1},\,
\alpha_0:=\delta+\epsilon_N-\epsilon_1~|~i=1,2,\ldots N-1\}~.
\end{equation}

Let $\mathfrak{sw}_{N+1}^{}$ be a subalgebra of $\mathfrak{gl}_{N+1}^{}$
generated by the root vectors: $e_{21}^{}$, $e_{i,i+1}^{}$, $e_{i+1,i}^{}$ for
$i=2,3,\ldots,N$ and $e_{N,N+1}^{}$, and also by the Cartan elements:
$e_{11}^{}+e_{22}^{}$, $e_{ii}^{}$ for $i=2,3,\ldots,N$. It is easy to check
that $\mathfrak{sw}_{N+1}$ has structure of a seaweed Lie algebra
(see \cite{DK}).

Let $\hat{\mathfrak{sw}}_{N}^{}$ be a subalgebra of $\mathfrak{gl}_{N}^{}[u]$
generated by the root vectors: $e_{21}^{(0)}$, $e_{i,i+1}^{(0)}$,
$e_{i+1,i}^{(0)}$ for $i=2,3,\ldots,N$ and $e_{N,1}^{(1)}$, and also by the
Cartan elements: $e_{ii}^{(0)}$ for $i=1,2,3,\ldots,N$. It is easy to check
that the Lie algebras $\hat{\mathfrak{sw}}_{N}$ and $\mathfrak{sw}_{N+1}^{}$
are isomorphic. This isomorphism is described by the following correspondence:
$e_{i+1,i}^{}\leftrightarrow e_{i+1,i}^{(0)}$ for $i=1,2,\ldots,N-1$,
$e_{i,i+1}^{}\leftrightarrow e_{i,i+1}^{(0)}$ for $i=2,3,\ldots,N-1$,
$e_{N,N+1}^{}\leftrightarrow e_{N,1}^{(1)}$ for $i=2,3,\ldots,N-1$, and
$(e_{11}^{}+e_{N+1,N+1}^{})\leftrightarrow e_{11}^{(0)}$,
$e_{ii}^{}\leftrightarrow e_{ii}^{(0)}$ for $i=2,3,\ldots,N$. We will call
$\hat{\mathfrak{sw}}_{N}$ an affine realization of $\mathfrak{sw}_{N+1}^{}$.

Now let us consider $q$-analogs of the previous Lie algebras. The quantum
algebra $U_{q}(\mathfrak{gl}_{N}^{})$ is generated by the Chevalley
elements\footnote{We denote the generators in the classical and quantum cases
by the same letter "$e$". It should not cause any misunderstanding.}
$e_{i,i+1}^{}$, $e_{i+1,i}^{}$ $(i=1,2,\ldots,N-1)$, $q^{\pm e_{ii}}$
$(i=1,2,\ldots,N)$ with the defining relations:
\begin{equation}\label{pr4}
\begin{array}{rcl}
q^{e_{ii}}q^{-e_{ii}}\!\!&=\!\!&q^{-e_{ii}}q^{e_{ii}}=1~,
\\[5pt]
q^{e_{ii}}q^{e_{jj}}\!\!&=\!\!&q^{e_{jj}}q^{e_{ii}}~,
\\[5pt]
q^{e_{ii}}e_{jk}^{}q^{-e_{ii}}\!\!&=\!\!&q^{\delta_{ij}-\delta_{ik}}
e_{jk}^{}\quad(|j-k|=1)~,
\\[3pt]
[e_{i,i+1}^{},\,e_{j+1,j}^{}]\!\!&=\!\!&\delta_{ij}\,
\mbox{\large$\frac{q^{e_{ii}-e_{i+1,i+1}}\,-\,q^{e_{i+1,i+1}-e_{ii}}}
{q\,-\,q^{-1}}$}~,
\\[5pt]
[e_{i,i+1}^{},\,e_{j,j+1}^{}]\!\!&=\!\!&0\quad{\rm for}\;\;|i-j|\geq 2~,
\\[5pt]
[e_{i+1,i}^{},\,e_{j+1,j}^{}]\!\!&=\!\!&0\quad{\rm for}\;\;|i-j|\geq 2~,
\\[5pt]
[[e_{i,i+1}^{},\,e_{j,j+1}^{}]_{q}^{},\,e_{j,j+1}^{}]_{q}^{}\!\!&=\!\!&0
\quad{\rm for}\;\;|i-j|=1~,
\\[5pt]
[[e_{i+1,i}^{},\,e_{j+1,j}^{}]_{q}^{},\,e_{j+1,j}^{}]_{q}^{}\!\!&=\!\!&0
\quad{\rm for}\;\;|i-j|=1~.
\end{array}
\end{equation}
where $[e_{\beta}^{},\,e_{\gamma}^{}]_{q}^{}$ denotes the $q$-commutator:
\begin{equation}\label{pr5}
[e_{\beta}^{},\,e_{\gamma}^{}]_{q}^{}\,:=\,e_{\beta}^{}e_{\gamma}^{}-
q^{(\beta,\gamma)}e_{\gamma}^{}e_{\beta}^{}~.
\end{equation}
The Hopf structure on $U_{q}(\mathfrak{gl}_{N}^{})$ is given by the following
formulas for comultiplication $\Delta_{q}$, antipode $S_{q}$, and co-unit
$\varepsilon_{q}$:
\begin{equation}\label{pr6}
\begin{array}{rcl}
\Delta_{q}(q^{\pm e_{ii}})\!\!&=\!\!&
q^{\pm e_{ii}}\otimes q^{\pm e_{ii}} ~,
\\[5pt]
\Delta_{q}(e_{i,i+1}^{})\!\!&=\!\!&
e_{i,i+1}^{}\otimes 1+q^{e_{i+1,i+1}-e_{ii}}\otimes e_{i,i+1}^{}~,
\\[5pt]
\Delta_{q}(e_{i+1,i}^{})\!\!&=\!\!&
e_{i+1,i}^{}\otimes q^{e_{ii}-e_{i+1,i+1}}+1\otimes e_{i+1,i}^{}~;
\end{array}
\end{equation}
\begin{equation}\label{pr7}
\begin{array}{rcl}
S_{q}(q^{\pm e_{ii}})\!\!&=\!\!&q^{\mp e_{ii}}~,
\\[5pt]
S_{q}(e_{i,i+1}^{})\!\!&=\!\!&-q^{e_{ii}-e_{i+1,i+1}}\,e_{i,i+1}^{}~,
\\[5pt]
S_{q}(e_{i+1,i}^{})\!\!&=\!\!&-e_{i+1,i}^{}\,q^{e_{i+1,i+1}-e_{i,i}}~;
\phantom{aaaaaaaaaa}
\end{array}
\end{equation}
\begin{equation}\label{pr8}
\begin{array}{rcccl}
\phantom{aa}\varepsilon_{q}(q^{\pm e_{ii}})^{}\!\!&=\!\!&1~,\quad
\varepsilon_{q}(e_{ij}^{}) \!\!&=\!\!&0\quad {\rm for}\;\;|i-j|=1~.
\end{array}
\end{equation}
For construction of the composite root vectors $e_{ij}^{}$ for $|i-j|\geq2$
we fix the following normal ordering of the positive root system $\Delta_{+}^{}$
(see \cite{T1,KT1,KT2})
\begin{equation}\label{pr9}
\begin{array}{c}
\epsilon_1^{}\!-\epsilon_2^{}\prec\epsilon_1^{}\!-\epsilon_3^{}
\prec\epsilon_2^{}\!-\epsilon_3^{}\prec\epsilon_{1}^{}\!-\epsilon_{4}^{}\prec
\epsilon_{2}^{}-\epsilon_{4}^{}\prec\epsilon_{3}^{}\!-\epsilon_{4}^{}
\prec\ldots\prec
\\[5pt]
\epsilon_{1}^{}\!-\epsilon_{k}^{}\!\prec\epsilon_{2}^{}\!-\epsilon_{k}^{}
\!\prec\ldots\prec\epsilon_{k-1}^{}\!-\epsilon_{k}^{}\!\prec\ldots\prec
\epsilon_{1}^{}\!-\epsilon_{N}^{}\!\prec\epsilon_{2}^{}\!-\epsilon_{N}^{}
\!\prec\ldots\prec\epsilon_{N-1}^{}\!-\epsilon_{N}^{}~.
\end{array}
\end{equation}
According to this ordering we set
\begin{equation}\label{pr10}
e_{ij}^{}\,:=\,[e_{ik}^{},\,e_{kj}^{}]_{q^{-1}}^{},\qquad
e_{ji}^{}\,:=\,[e_{jk}^{},\,e_{ki}^{}]_{q}^{}~,
\end{equation}
where $1\le i<k<j\le N$. It should be stressed that the structure of the
composite root vectors is not independent on a choice of the index $k$ in the
r.h.s. of the definition (\ref{pr10}). In particular, we have
\begin{equation}\label{pr11}
\begin{array}{rcccl}
e_{ij}^{}\!\!&:=\!\!\!&[e_{i,i+1}^{},\,e_{i+1,j}^{}]_{q^{-1}}^{}&\!\!\!=
\!\!\!&[e_{i,j-1}^{},\,e_{j-1,j}^{}]_{q^{-1}}^{}~,
\\[7pt]
e_{ji}^{}&\!\!:=\!\!&[\,e_{j,i+1}^{},\,e_{i+1,i}^{}\,]_{q}^{}&\!\!=
\!\!&[e_{j,j-1}^{},\,e_{j-1,i}^{}]_{q}^{}~,
\end{array}
\end{equation}
where $2\le i+1<j\le N$.

Using these explicit constructions and the defining relations (\ref{pr4}) for
the Chevalley basis it is not hard to calculate the following relations between
the Cartan--Weyl generators $e_{ij}$ ($i,j=1,2,\ldots, N$): 
\begin{eqnarray}
q^{e_{kk}^{}}e_{ij}^{}q^{-e_{kk}^{}}\!\!&=\!\!&q^{\delta_{ki}^{}-\delta_{kj}^{}}
e_{ij}^{}\qquad(1\le i,j,k\le N)~,
\label{pr12}
\\[5pt]
[e_{ij}^{},\,e_{ji}^{}]\!\!&=\!\!&\frac{q^{e_{ii}^{}-e_{jj}^{}}-
q^{e_{jj}^{}-e_{ii}^{}}}{q-q^{-1}}\qquad(1\le i<j\le N)~,
\label{pr13}
\\[5pt]
[e_{ij}^{},\,e_{kl}^{}]_{q^{-1}}\!\!&=\!\!&\delta_{jk}^{}e_{il}^{}
\qquad(1\le i<j\le k<l\le N)~,
\label{pr14}
\\[5pt]
[e_{ik}^{},\,e_{jl}^{}]_{q^{-1}}^{}\!\!&=\!\!&(q-q^{-1})\,e_{jk}^{}e_{il}^{}
\qquad(1\le i<j<k<l\le N)~,
\label{pr15}
\\[5pt]
[e_{jk}^{},\,e_{il}^{}]_{q^{-1}}^{}\!\!&=\!\!&0\qquad(1\le i\le j<k\le l\le N)~,
\label{pr16}
\\[5pt]
[e_{kl}^{},\,e_{ji}^{}]\!\!&=\!\!&0\qquad(1\le i<j\le k<l\le N)~,
\label{pr17}
\\[5pt]
[e_{il}^{},\,e_{kj}^{}]\!\!&=\!\!&0\qquad(1\le i<j<k<l\le N)~,
\label{pr18}
\\[5pt]
[e_{ji}^{},\,e_{il}^{}]\!\!&=\!\!&e_{jl}^{}\,q^{e_{ii}^{}-e_{jj}^{}}
\qquad(1\le i<j<l\le N)~,
\label{pr19}
\\[5pt]
[e_{kl}^{},\,e_{li}^{}]\!\!&=\!\!&e_{ki}^{}\,q^{e_{kk^{}}-e_{ll}^{}}
\qquad(1\le i<k<l\le N)~,
\label{pr20}
\\[5pt]
[e_{jl}^{},\,e_{ki}^{}]\!\!&=\!\!&(q^{-1}-q)\,e_{kl}^{}e_{ji}^{}\,
q^{e_{jj}-e_{kk}}\qquad(1\le i<j<k<l\le N)~.
\label{pr21}
\end{eqnarray}
These formulas can be obtained from the relations between elements of the Cartan--Weyl
basis for  the quantum superalgebra $U_q(\mathfrak{gl}(N|M)$ (see \cite{T3}).
If we apply the Cartan involution ($e_{ij}^{*}=e_{ji}^{}$) the formulas above,
we will get all relations between elements of the Cartan--Weyl basis.

The quantum algebra $U_{q}(\mathfrak{gl}_{N}^{}[u])$ ($N\ge3$) is generated (as
an unital associative algebra over $\mathbb C[[\log q]]$) by the algebra
$U_{q}(\mathfrak{gl}_{N}^{})$ and the additional element $e_{N1}^{(1)}$ with
the relations:
\begin{equation}\label{pr22}
\begin{array}{rcl}
q^{\pm e_{ii}^{(0)}}e_{N1}^{(1)}\!\!&=\!\!&q^{\mp(\delta_{i1}-\delta_{iN})}
e_{N1}^{(1)}q^{\pm e_{ii}^{(0)}}~,
\\[5pt]
[e_{i,i+1}^{(0)},~e_{N1}^{(1)}]\!\!&=\!\!&0 \quad{\rm for}\;\, i=2,3,\ldots,N-2~,
\\[5pt]
[e_{i+1,i}^{(0)},~e_{N1}^{(1)}]\!\!&=\!\!&0\quad{\rm for}\;\, i=1,2,\ldots,N-1~,
\\[5pt]
[e_{12}^{(0)},~[e_{12}^{(0)},~e_{N1}^{(1)}]_{q}^{}]_{q}^{}\!\!&=\!\!&0~,
\\[5pt]
[e_{N-1,N}^{(0)},~[e_{N-1,N}^{(0)},~e_{N1}^{(1)}]_{q}^{}]_{q}^{}\!\!&=\!\!&0~,
\\[5pt]
[[e_{12}^{(0)},~e_{N1}^{(1)}]_{q}^{},~e_{N1}^{(1)}]_{q}^{}\!\!&=\!\!&0~,
\\[5pt]
[[e_{N-1,N}^{(0)},~e_{N1}^{(1)}]_{q}^{},~e_{N1}^{(1)}]_{q}^{}\!\!&=\!\!&0~,
\end{array}
\end{equation}
The Hopf structure of  $U_{q}(\mathfrak{gl}_{N}^{}[u])$ is defined by the
formulas (\ref{pr6})-(\ref{pr8}) for $U_{q}(\mathfrak{gl}_{N}^{(0)})$ and the
following additional formulas for the comultiplication and the antipode:
\begin{eqnarray}\label{pr23}
\Delta_{q}(e_{N1}^{(1)})\!\!&=\!\!&e_{N1}^{(1)}\otimes 1+
q^{e_{11}^{(0)}-e_{NN}^{(0)}}\otimes e_{N1}^{(1)}~,
\\[5pt]\label{pr24}
S_{q}(e_{N1}^{(1)})\!\!&=\!\!&-q^{e_{NN}^{(0)}-e_{11}^{(0)}}e_{N1}^{(1)}~.
\end{eqnarray}

Quantum analogs of the seaweed algebra $\mathfrak{sw}_{N+1}^{}$ and its affine
realization $\hat\mathfrak{sw}_{N}^{}$ are inherited from the quantum algebras
$U_{q}(\mathfrak{gl}_{N+1}^{})$ and $U_{q}(\mathfrak{gl}_{N}^{}[u])$. Namely,
the quantum algebra $U_q(\mathfrak{sw}_{N+1}^{})$ is generated by the root
vectors: $e_{21}^{}$, $e_{i,i+1}^{}$, $e_{i+1,i}^{}$ for $i=2,3,\ldots,N$ and
$e_{N,N+1}^{}$, and also by the $q$-Cartan elements:
$q^{e_{11}^{}+e_{N+1,N+1}^{}}$, $q^{e_{ii}^{}}$ for $i=2,3,\ldots,N$ with the
relations satisfying (\ref{pr4}). Similarly,  the quantum algebra
$U_q(\hat\mathfrak{sw}_{N}^{})$ is generated by the root vectors:
$e_{21}^{(0)}$, $e_{i,i+1}^{(0)}$, $e_{i+1,i}^{(0)}$ for $i=2,3,\ldots,N$ and
$e_{N,1}^{(1)}$, and also by the $q$-Cartan elements: $q^{e_{ii}^{(0)}}$ for
$i=1,2,3,\ldots,N$ with the relations satisfying (\ref{pr4}) and (\ref{pr22}).
It is clear that the algebras $U_{q}(\mathfrak{gl}_{N+1}^{})$ and
$U_q(\hat\mathfrak{sw}_{N}^{})$ are isomorphic as associative algebras but they
are not isomorphic as Hopf algebras. However if we introduced a new coproduct
in the Hopf algebra $U_{q}(\mathfrak{gl}_{N+1}^{})$
\begin{equation}\label{pr25}
\Delta_{\;q}^{(\mathfrak{F}_{1,N+1}^{})}(x)\,=\,\mathfrak{F}_{1,N+1}^{}
\Delta_q(x)\mathfrak{F}_{1,N+1}^{-1}\quad (\forall x\in
U_{q}(\mathfrak{gl}_{N+1}^{}))~,
\end{equation}
where
\begin{equation}\label{pr26}
\mathfrak{F}_{1,N+1}^{}:=q^{\,-e_{11}^{}\otimes e_{N+1,N+1}^{}}~,
\end{equation}
we obtain an isomorphism of Hopf algebras
\begin{equation}\label{pr27}
U_q^{(\mathfrak{F}_{1,N+1})}(\mathfrak{sw}_{N+1}^{})\,\simeq\,
U_q(\hat\mathfrak{sw}_{N}^{})~.
\end{equation}
Here the symbol $U_q^{(\mathfrak{F}_{1,N+1})}(\mathfrak{sw}_{N+1}^{})$ denotes
the quantum seaweed algebra $U_q(\mathfrak{sw}_{N+1}^{})$ with the twisted
coproduct (\ref{pr25}).

\section{Cartan part of Cremmer-Gervais $r$-matrix}

First of all we recall classification of quasi-triangular $r$-matrices for
a simple Lie algebra $\mathfrak{g}$. The quasi-triangular $r$-matrices are
solutions of the system
\begin{equation}\label{cg1}
\begin{array}{rcl}
r^{12}+r^{21}\!\!&=\!\!&\Omega~,
\\[7pt]
[r^{12},r^{13}]+[r^{12},r^{23}]+[r^{13},r^{23}]\!\!&=\!\!&0~,
\end{array}
\end{equation}
where $\Omega$ is the quadratic the Casimir two-tensor in
$\mathfrak{g}\otimes\mathfrak{g}$. Belavin and Drinfeld proved that any
solution of this system is defined by a triple $(\Gamma_{1},\Gamma_{2}, \tau)$,
where $\Gamma_{1},\Gamma_{2}$ are subdiagrams of the Dynkin diagram of
$\mathfrak{g}$ and $\tau$ is an isometry between these two subdiagrams.
Further, each $\Gamma_i$ defines a reductive subalgebra of $\mathfrak{g}$, and
$\tau$ is extended to an isometry (with respect to the corresponding
restrictions of the Killing form) between the corresponding reductive
subalgebras of $\mathfrak{g}$. The following property of $\tau$ should be
satisfied: $\tau^{k}(\alpha)\not\in\Gamma_{1}$ for any $\alpha\in\Gamma_{1}$
and some $k$. Let $\Omega_0$ be the Cartan part of $\Omega$. Then one can
construct a quasi-triangular $r$-matrix according to the following
\begin{theorem}[Belavin--Drinfeld \cite{BD}]
Let $r_{0}\in\mathfrak{h}\otimes\mathfrak{h}$ satisfies the systems
\begin{eqnarray}\label{cg2}
r_{0}^{12}+r_{0}^{21}\!\!&=\!\!&\Omega_{0}~,
\\[5pt]\label{cg3}
(\alpha\otimes 1+1\otimes\alpha)(r_{0})\!\!&=\!\!&h_{\alpha}^{}
\\[5pt]\label{cg4}
(\tau(\alpha)\otimes 1+1\otimes\alpha)(r_{0})\!\!&=\!\!&0
\end{eqnarray}
for any $\alpha\in\Gamma_{1}$. Then the tensor
\begin{eqnarray}\label{cg5}
r=r_{0}+\sum_{\alpha>0}e_{-\alpha}^{}\otimes e_{\alpha}^{}+
\sum_{\alpha>0;k\geq1}e_{-\alpha}^{} \wedge e_{\tau^k(\alpha)}
\end{eqnarray}
satisfies (\ref{cg1}). Moreover, any solution of the system (\ref{cg1}) is of
the above form, for a suitable triangular decomposition of $\mathfrak{g}$ and
suitable choice of a basis $\{e_{\alpha}^{}\}$.
\end{theorem}
In what follows, for aim of quantization of algebra structures on the
polynomial Lie algebra $\mathfrak{gl}_{N}^{}[u])$ we will use the twisted
two-tensor $q^{r_{0}^{}(N)}$ where $r_{0}^{}(N)$ is the Cartan part of the
Cremmer--Gervais $r$-matrix for the Lie algebra $\mathfrak{gl}_{N}^{}$ when
$\Gamma_{1}=\{\alpha_{1},\alpha_{2},\ldots,\alpha_{N-2}\}$
$\Gamma_{2}=\{\alpha_{2},\alpha_{3},\ldots,\alpha_{N-1}\}$ and
$\tau(\alpha_i)=\alpha_{i+1}$. An explicit form of $r_{0}^{}(N)$ is defined by
the following proposition (see \cite{GG}).
\begin{proposition}
The Cartan part of the Cremmer--Gervais $r$-matrix for $\mathfrak{gl}_{N}^{}$
is given by the following expression
\begin{equation}\label{cg6}
r_{0}^{}(\mathfrak{gl}_{N}^{})\,=\,
\displaystyle\frac{1}{2}~\sum_{i=1}^{N}e_{ii}^{}\otimes e_{ii}^{}+
\sum_{1\leq i<j\leq N}\,\frac{N+2(i-j)}{2N}\;e_{ii}^{}\wedge e_{jj}^{}~.
\end{equation}
\end{proposition}
It is easy to check that the Cartan part (\ref{cg6}), $r_{0}^{}(N):=
r_{0}^{}(\mathfrak{gl}_{N}^{})$, satisfies the conditions
\begin{eqnarray}\label{cg7}
\bigl(\epsilon_k^{}\otimes{\rm id}+{\rm id}\otimes\epsilon_k^{}\bigr)
\bigl(r_{0}^{}(N)\bigr)\!\!&=\!\!&e_{kk}^{}\quad\;\;{\rm for}\;\;k=1,2,\ldots,N,
\\[5pt]\label{cg8}
\bigl(\epsilon_{k}^{}\otimes{\rm id}+{\rm id}\otimes\epsilon_{k'}^{}\bigr)
\bigl(r_{0}^{}(N)\bigr)\!\!&=\!\!&(k-k')\,\mathcal{C}_1^{}(N)-\!\!\!
\sum_{i=k'+1}^{k-1}e_{ii}^{}\quad\;{\rm for}\;\; 1\leq k'<k\leq N,
\end{eqnarray}
where $\mathcal{C}_1^{}(N)$ is the normalized central element:
\begin{equation}\label{cg9}
\mathcal{C}_1^{}(N)\,:=\,\frac{1}{N}\sum_{i=1}^{N} e_{ii}^{}~.
\end{equation}
In particular (\ref{cg7}) and (\ref{cg8}) imply the Belavin--Drinfeld
conditions (\ref{cg3}) and (\ref{cg4}), i.e.
\begin{eqnarray}\label{cg10}
\bigl(\alpha_k^{}\otimes{\rm id}+{\rm id}\otimes\alpha_k^{}\bigr)
\bigl(r_{0}^{}(N)\bigr)\!\!&=\!\!& h_{\alpha_k}^{}\,:=\,e_{kk}^{}-
e_{k+1,k+1}^{}~,
\\[5pt]\label{cg11}
\bigl(\tau(\alpha_{k'}^{})\otimes{\rm id}+{\rm id}\otimes\alpha_{k'}^{}\bigr)
\bigl(r_{0}^{}(N)\bigr)\!\!&=\!\!&\bigl(\alpha_{k'+1}^{}\otimes{\rm id}+ {\rm
id}\otimes\alpha_{k'}^{}\bigr)\bigl(r_{0}^{}(N)\bigr)\;=\;0
\end{eqnarray}
for $k=1,2\ldots,N-1~$, $k'=1,2\ldots,N-2~$, where $\alpha_k^{}=\epsilon_k^{}-
\epsilon_{k+1}^{}$ and $\alpha_{k'}^{}=\epsilon_{k'}^{}-\epsilon_{k'+1}^{}$ are
the simple roots of system $\mathop{\Pi}(\mathfrak{gl}_{N})$ (see (\ref{pr3})).

Now we consider some properties of the two-tensor $q^{r_{0}^{}(N)}$. First of
all it is evident that this two-tensor satisfies cocycle equation (\ref{i4}).
Further, for construction of a twisting two-tensor corresponding to the
Cremmer--Gervais $r$-matrix (\ref{cg6}) we introduce new Cartan--Weyl basis
elements $e_{ij}^{\,\prime}$ ($i\neq j$) for the quantum algebra
$U_{q}(\mathfrak{gl}_{N})$ as follows
\begin{eqnarray}\label{cg12}
e_{ij}^{\,\prime}\!\!&=\!\!&\displaystyle e_{ij}^{}q^{\,\bigl((\epsilon_i-
\epsilon_{j})\otimes {\rm id}\bigr)\bigl(r_{0}(N)\bigr)}\,=
\,e_{ij}^{}\,q^{\;\sum\limits_{k=i}^{j-1}e_{kk}^{}-(j-i)\mathcal{C}_{1}^{}(N)}~,
\\[5pt]\label{cg13}
e_{ji}^{\,\prime}\!\!&=\!\!&q^{\,\bigl({\rm id}\otimes(\epsilon_{j}-
\epsilon_{i})\bigr)\bigl(r_{0}(N)\bigr)}\,e_{ji}^{}\,=\,
q^{\;\sum\limits_{k=i+1}^{j}e_{kk}^{}-(j-i)\mathcal{C}_{1}^{}(N)}\,e_{ji}^{}~,
\end{eqnarray}
for $1\leq i<j\leq N$. Permutation relations for these elements can be 
obtained from the relations (\ref{pr12})--(\ref{pr21}). For example, we have
\begin{equation}\label{cg14}
\begin{array}{rcl}
[e_{ij}^{\prime},\,e_{ji}^{\prime}]\!\!&=\!\!&\displaystyle
[e_{ij}^{},\,e_{ji}^{}]\,q^{\,\bigl((\epsilon_i-\epsilon_{j})\otimes{\rm
id}+{\rm id}\otimes(\epsilon_{j}-\epsilon_{i})\bigr)\bigl(r_{0}^{}(N)\bigr)}
\\[12pt]
\!\!&=\!\!&\displaystyle\frac{q^{\,2\!\sum\limits_{k=i}^{j-1}e_{kk}^{}-
2(j-i)\mathcal{C}_{1}^{}(N)}-q^{\,2\!\!\!\sum\limits_{k=i+1}^{j}e_{kk}^{}-
2(j-i)\mathcal{C}_{1}^{}(N)}}{q-q^{-1}}~.
\end{array}
\end{equation}
It is not hard to check that the Chevalley elements $e_{i,i+1}^{\prime}$ and
$e_{i+1,i}^{\prime}$ have the following coproducts after twisting by the
two-tensor $q^{r_{0}^{}(N)}$:
\begin{equation}\label{cg15}
\begin{array}{rcl}
q^{\,r_{0}^{}(N)}\Delta_q^{}(e_{i,i+1}^{\prime})q^{-r_{0}^{}(N)}\!\!&=\!\!&
e_{i,i+1}^{\prime}\otimes q^{\,2\bigl((\epsilon_i-\epsilon_{i+1})
\otimes{\rm id}\bigr)(r_{0}^{}(N))}+1\otimes e_{i,i+1}^{\prime}
\\[5pt]
\!\!&=\!\!&e_{i,i+1}^{\prime}
\otimes q^{\,2e_{ii}^{}-2\mathcal{C}_{1}^{}(N)}+1\otimes e_{i,i+1}^{\prime}~,
\end{array}
\end{equation}
\begin{equation}\label{cg16}
\begin{array}{rcl}
q^{\,r_{0}^{}(N)}\Delta_q^{}(e_{i+1,i}^{\prime})q^{-r_{0}^{}(N)}\!\!&=\!\!&
e_{i+1,i}^{\prime}\otimes1+q^{-2\bigl({\rm id}\otimes(\epsilon_{i+1}-
\epsilon_{i})\bigr)(r_{0}^{}(N))}\otimes e_{i+1,i}^{\prime}
\\[5pt]
\!\!&=\!\!&e_{i+1,i}^{\prime}\otimes1+
q^{2e_{i+1,i+1}^{}-2\mathcal{C}_{1}^{}(N)}\otimes e_{i+1,i}^{\prime}~.
\end{array}
\end{equation}
for $1\leq i<N$. Since the quantum algebra $U_{q}(\mathfrak{gl}_{N})$ is a
subalgebra of the quantum affine algebra $U_{q}(\mathfrak{gl}_{N}[u])$ let us
introduce the new affine root vector $e_{N1}^{\prime(1)}$ in accordance with
(\ref{cg12}):
\begin{eqnarray}\label{cg17}
e_{N1}^{\prime(1)}\!\!&=\!\!&\displaystyle e_{N1}^{}q^{\,\bigl((\epsilon_{N}^{}-
\epsilon_{1}^{})\otimes {\rm id}\bigr)\bigl(r_{0}(N)\bigr)}\,=
\,e_{N1}^{}\,q^{e_{NN}^{}-\mathcal{C}_{1}^{}(N)}~.
\end{eqnarray}
The coproduct of this element after twisting by the two-tensor $q^{r_{0}^{}(N)}$
has the form
\begin{equation}\label{cg18}
\begin{array}{rcl}
q^{\,r_{0}^{}(N)}\Delta_q^{}(e_{N1}^{\prime(1)})q^{-r_{0}^{}(N)}\!\!&=\!\!&
e_{N1}^{\prime(1)}\otimes q^{\,2\bigl((\epsilon_{1}^{}-\epsilon_{N}^{})
\otimes{\rm id}\bigr)(r_{0}^{}(N))}+1\otimes e_{N1}^{\prime(1)}
\\[7pt]
\!\!&=\!\!&e_{N1}^{\prime(1)}\otimes q^{\,2e_{NN}^{}-2\mathcal{C}_{1}^{}(N)}+
1\otimes e_{N1}^{\prime(1)}~.
\end{array}
\end{equation}

Consider the quantum seaweed algebra $U_{q}^{}(\mathfrak{sw}_{N+1}^{})$ after
twisting by the two-tensor $q^{r_{0}^{}(N+1)}$. Its new Cartan--Weyl basis and
the coproduct for the Chevalley generators are given by formulas (\ref{cg12}),
(\ref{cg13}) and (\ref{cg15}), (\ref{cg16}), where $N$ should be replaced  by
$N+1$ and where $i\neq1$ in (\ref{cg12}) and (\ref{cg15}), and  $j\neq N$ in
(\ref{cg13}), and $i\neq N$ in (\ref{cg16}). In particular, for the element
$e_{N,N+1}^{\prime}$ we have
\begin{eqnarray}\label{cg19}
e_{N,N+1}^{\prime}\!\!&=\!\!&e_{N,N+1}^{}\,q^{e_{NN}^{}-
\mathcal{C}_{1}^{}(N+1)}~.
\\[5pt]
\label{cg20}
q^{\,r_{0}^{}(N+1)}\Delta_q^{}(e_{N,N+1}^{\prime})q^{-r_{0}^{}(N+1)}\!\!&=
\!\!&e_{N,N+1}^{\prime}\otimes q^{\,2e_{NN}^{}-2\mathcal{C}_{1}^{}(N+1)}+
1\otimes e_{N,N+1}^{\prime}~.
\end{eqnarray}
Comparing the Hopf structure of the quantum seaweed algebra
$U_{q}^{}(\mathfrak{sw}_{N+1}^{})$  after twisting by the two-tensor
$q^{r_{0}^{}(N+1)}$and its affine realization
$U_{q}^{}(\hat\mathfrak{sw}_{N}^{})$ after twisting by the two-tensor
$q^{r_{0}^{}(N)}$ we see that these algebras are isomorphic as Hopf
algebras:
\begin{eqnarray}\label{cg21}
q^{\,r_{0}^{}(N+1)}\Delta_q^{}(U_{q}^{}(\mathfrak{sw}_{N+1}^{}))
q^{-r_{0}^{}(N+1)}\!\!&\simeq\!\!&q^{\,r_{0}^{}(N)}\Delta_q^{}
(U_{q}^{}(\hat\mathfrak{sw}_{N}^{}))q^{-r_{0}^{}(N)}~.
\end{eqnarray}
In terms of new Cartan--Weyl bases this isomorphism, $"\imath"$, is arranged as
follows
\begin{eqnarray}\label{cg22}
\imath(e_{ij}^{\,\prime})\!\!&=\!\!&e_{ij}^{\prime(0)}
\qquad{\rm for}\;\;2\leq i<j\leq N~,
\\[7pt]\label{cg23}
\imath(e_{ji}^{\,\prime})\!\!&=\!\!&e_{ji}^{\prime(0)}
\qquad{\rm for}\;\;1\leq i<j\leq N-1~,
\\[7pt]\label{cg24}
\imath(e_{ii}^{}-\mathcal{C}_{1}^{}(N+1))\!\!&=\!\!&
e_{ii}^{(0)}-\mathcal{C}_{1}^{}(N)~,\qquad{\rm for}\;\;2\leq i\leq N~,
\\[2pt]\label{cg25}
\imath(e_{iN+1}^{\,\prime})\;=\;e_{i1}^{\prime(1)}\!\!&=\!\!&e_{i1}^{(1)}
q^{\,\bigl((\epsilon_{i}^{}-\epsilon_{1}^{})\otimes{\rm id}\bigr)
\bigl(r_{0}(N)\bigr)}\;=\;e_{i1}^{(1)}q^{\;\sum\limits_{k=i}^{N}
e_{kk}^{}-(N+1-i)\mathcal{C}_{1}^{}(N)}
\end{eqnarray}
for $2\leq i\leq N$. Where the affine root vectors $e_{i1}^{(1)}$
($2\leq i<N$) are defined by the formula (cf. \ref{pr10}):
\begin{equation}\label{cg26}
e_{i1}^{(1)}\,=\,[e_{iN}^{(0)},\,e_{N1}^{(1)}]_{q^{-1}}^{}~.
\end{equation}

\section{Affine realization of a Cremmer-Gervais twist}

For construction of a twisting two-tensor corresponding to the Cremmer-Gervais
$r$-matrix (\ref{cg6}) we will follow to the papers \cite{ESS, IO}.

Let $\mathcal{R}$ be a universal $R$-matrix of the quantum algebra
$U_{q}(\mathfrak{gl}_{N+1})$. According to \cite{KT1} it has the following
form
\begin{equation}\label{cgt1}
\mathcal{R}\,=\,R\cdot K
\end{equation}
where the factor $K$ is a $q$-power of Cartan elements (see \cite{KT1}) and we
do not need its explicit form. The factor $R$ depends on the root vectors and
it is given by the following formula
\begin{equation}\label{cgt2}
\begin{array}{rcl}
R\!\!&=\!\!&R_{12}^{}(R_{13}^{}R_{23}^{})(R_{14}^{}R_{24}^{}R_{34}^{})
\cdots(R_{1,N+1}^{}R_{2,N+1}^{}\cdots R_{N,N+1}^{})
\\[5pt]
\!\!&=\!\!&\displaystyle\uparrow\prod_{j=2}^{N+1}\Bigl(\uparrow
\prod_{i=1}^{j-1}R_{ij}^{}\Bigr)~,
\end{array}
\end{equation}
where
\begin{equation}\label{cgt3}
R_{ij}^{}\,=\,\exp_{q^{-2}}((q-q^{-1})e_{ij}^{}\otimes e_{ji}^{})~,
\end{equation}
\begin{equation}\label{cgt4}
\exp_{q}(x):=\sum_{n\ge 0}\frac{x^{n}}{(n)_{q}!}~,
\quad(n)_{q}!\equiv (1)_{q}(2)_{q}\ldots (n)_{q}~,\quad
(k)_{q}\equiv (1-q^{k})/(1-q)~.
\end{equation}
It should be noted that the product of factors $R_{ij}^{}$ in (\ref{cgt2})
corresponds to the normal ordering (\ref{pr9}) where $N$ is replaced by $N+1$.

Let $R^{\,\prime}:=q^{r_{0}^{}(N+1)}R\,q^{-r_{0}^{}(N+1)}$. It is evident that
\begin{equation}\label{cgt5}
R^{\,\prime}\,=\,\uparrow\!\!\prod_{j=2}^{N+1}\Bigl(\uparrow\!
\prod_{i=1}^{j-1}R_{ij}^{\,\prime}\Bigr)~,
\end{equation}
where
\begin{equation}\label{cgt6}
R_{ij}^{\,\prime}\,=\,\exp_{q^{-2}}((q-q^{-1})e_{ij}^{\,\prime}\otimes
e_{ji}^{\,\prime})~.
\end{equation}
Here $e_{ij}^{\,\prime}$ and $e_{ji}^{\,\prime}$ are the root vectors
(\ref{cg12}) and (\ref{cg13}) where $N$ should be replace by $N+1$.

Let $\mathcal{T}$ be a homomorphism operator which acts on the elements
$e_{ij}^{\,\prime}$ ($1\leq i<j\leq N+1$) by formulas
$\mathcal{T}(e_{ij}^{\,\prime})=e_{\tau(ij)}^{\,\prime}=e_{i+1,j+1}^{\,\prime}$
for $1\leq i<j\leq N$, and $\mathcal{T}(e_{i,N+1}^{\,\prime})=0$ for all
$i=1,2,\ldots, N$. We set
\begin{equation}\label{cgt7}
R^{\,\prime(k)}\,:=\,(\mathcal{T}^k\otimes{\rm id})(R^{\,\prime})\,=\,\uparrow
\!\!\!\prod_{j=2}^{N+1-k}\Bigl(\uparrow\prod_{i=1}^{j-1}R_{ij}^{\,\prime(k)}\Bigr)~,
\end{equation}
where
\begin{equation}\label{cgt8}
R_{ij}^{\,\prime(k)}\,=\,\exp_{q^{-2}}\bigl((q-q^{-1})
\mathcal{T}^k(e_{ij}^{\,\prime})\otimes e_{ji}^{\,\prime}\bigr)\,=\,
\exp_{q^{-2}}\bigl((q-q^{-1})e_{i+k,j+k}^{\,\prime}\otimes e_{ji}^{\,\prime}\bigr)
\end{equation}
for $k\leq N-j$.

According to \cite{ESS, IO}, the Cremmer-Gervais twist $\mathcal{F}_{CG}^{}$ in
$U_{q}(\mathfrak{gl}_{N+1})$ is given as follows
\begin{equation}\label{cgt9}
\mathcal{F}_{CG}^{}\,=\,F\cdot q^{r_{0}^{}(N+1)}~,
\end{equation}
where
\begin{equation}\label{cgt10}
F\,=\,R^{\,\prime(N-1)}R^{\,\prime(N-2)}\cdots R^{\,\prime(1)}~.
\end{equation}
It is easy to see that the support of the twisting two-tensor (\ref{cgt10}) is
the quantum seaweed algebra $U_{q}^{}(\mathfrak{sw}_{N+1}^{})$ with the
coproducts (\ref{cg15}) and (\ref{cg16}) where $N$ should be replace by $N+1$.
From the results of the previous section it follows that we can immediately
obtain an affine realization $\hat\mathcal{F}_{CG}^{}$ which twists the quantum
affine algebra $U_{q}^{}(\mathfrak{gl}_{N}^{}[u])$:
\begin{eqnarray}\label{cgt11}
\hat\mathcal{F}_{CG}^{}\!\!&=\!\!&\hat{F}\cdot q^{r_{0}^{}(N)}~,
\\[7pt]\label{cgt12}
\hat{F}\,:=\,(\imath\otimes\imath)(F)\!\!&=\!\!&\hat{R}^{\,\prime(N-1)}
\hat{R}^{\,\prime(N-2)}\cdots\hat{R}^{\,\prime(1)}
\\[5pt]\label{cgt13}
\hat{R}^{\,\prime(k)}\!\!&=\!\!&\uparrow\!\!\!\prod_{j=2}^{N+1-k}
\Bigl(\uparrow\prod_{i=1}^{j-1}\hat{R}_{ij}^{\,\prime(k)}\Bigr)~.
\end{eqnarray}
where
\begin{eqnarray}\label{cgt14}
\hat{R}_{ij}^{\,\prime(k)}\!\!&=\!\!&\exp_{q^{-2}}\Bigl((q-q^{-1})
e_{i+k,j+k}^{\,\prime(0)}\otimes e_{ji}^{\,\prime(0)}\Bigr)\qquad
{\rm for}\;\;1\leq i<j\leq N-k~,
\\[5pt]\label{cgt15}
\hat{R}_{i,N+1-k}^{\,\prime(k)}\!\!&=\!\!&\exp_{q^{-2}}
\Bigl((q-q^{-1})e_{N,i+k}^{\,\prime(1)}\otimes
e_{N+1-k,i}^{\,\prime(0)}\Bigr)\quad{\rm for}\;\;1\leq i<N+1-k~.
\end{eqnarray}

\section{From constant to affine twists}

\subsection{Affine $r$-matrices of Cremmer-Gervais type and their quantization}

In this subsection we are going to describe two classes of the
quasi-trigonometric $r$-matrices of Cremmer-Gervais type, which we call affine
$r$-matrices of the Cremmer-Gervais type.

We need some results proved in \cite{KPST} for the case $\mathfrak{sl}_N$. Let
us consider the Lie algebra $\mathfrak{sl}_N\oplus \mathfrak{sl}_N$ with the
form $Q((a,b),(c,d))=K(a,c)-K(b,d)$, where $K$ is the Killing form on
$\mathfrak{sl}_N$. Let $P_i$ (resp. $P_i^-$ ) be the maximal parabolic
subalgebra containing all positive (resp. negative) roots and not containing
all negative roots which have the simple root $\alpha_i$ in their
decomposition. Assume we have two parabolic (or maybe $\mathfrak{sl}_N)$
subalgebras $S_1,\ S_2$ of $\mathfrak{sl}_N$ such that their reductive parts
$R_1,\ R_2$ are isomorphic and let $T: R_1\to R_2$ be an isomorphism, which is
an isometry with respect to the reductions of the Killing form $K$ onto $R_1$
and $R_2$. Then the triple $(S_1,S_2,T)$ defines a Lagrangian subalgebra of
$\mathfrak{sl}_N\oplus \mathfrak{sl}_N$ (with respect to the form $Q$). It was
proved in \cite{KPST} that any quasi-trigonometric solution is defined by a
Lagrangian subalgebra $W\subset  \mathfrak{sl}_N\oplus \mathfrak{sl}_N  $ such
that $W$ is transversal to $(P_i,P_i,\rm{id})$. Now we would like to present
two types of $W$, which define quasi-trigonometric $r$-matrices for
$\mathfrak{sl}_N$ with $i=N-1$.

Type 1, which we call {\it affine Cremmer--Gervais I classical $r$-matrix}: In
this case $S_1=S_2=\mathfrak{sl}_N$, $T=\rm{Ad}(\sigma_N)$, where $\sigma_N$ is
the permutation matrix, which represents the cycle
$\{1,2,\ldots,N\}\to\{2,3,\ldots,N,1\}$.
\begin{proposition}\label{p2}
The triple $(S_1,S_2,T)$ above defines a quasi-trigonometric $r$-matrix.
\end{proposition}
Type 2, which we call {\it affine Cremmer--Gervais II classical $r$-matrix}: In
this case $S_3=P_{1}^- ,\ S_4=P_{N-1}^-,\ T'(e_{ij})=e_{i-1,j-1}$, where
$e_{ij}$ are the matrix units.
\begin{proposition}\label{p3}
The triple $(S_3,S_4,T')$ above defines a quasi-trigonometric $r$-matrix.
\end{proposition}
Both propositions can be proved directly. 
\begin{theorem}\label{th2}Let $\hat{{\cal F}}_{CG}^{'}$ be the
twist (\ref{cgt9}) reduced to $U_{q}(\hat{\mathfrak{sl}}_{N})$, and let
$\hat{{\cal R}}$ be the universal $R$-matrix for
$U_{q}(\hat{\mathfrak{sl}}_{N})$. Then the $R$-matrix $\hat{{\cal
F}}_{CG}^{'21}\hat{{\cal R}}\hat{{\cal F}}_{CG}^{'-1}$ quantizes the
quasi-trigonometric $r$-matrix obtained in Proposition \ref{p2}.
\end{theorem}
Now we turn to quantization of the quasi-trigonometric $r$-matrix from
Proposition \ref{p3}. In order to do this we propose another method for
construction of affine twists from constant ones. Let
$\{\alpha_{0},\alpha_{1},\ldots,\alpha_{N-1}\}$ be the vertices of the Dynkin
diagram of $\mathfrak{sl}_{N}$. Then the map $\tau$:
$\{\alpha_{1},\ldots,\alpha_{N-1}\}$ to
$\{\alpha_{2},\ldots,\alpha_{N-1},\alpha_{0}\}$ defines an embedding of Hopf
algebras $U_{q}(\mathfrak{sl}_{N})\hookrightarrow
U_{q}(\hat{\mathfrak{sl}}_{N})$. Abusing notations we denote this embedding by
$\tau$. The map $\tau$ sends each twist ${\cal F}$ of
$U_{q}(\mathfrak{sl}_{N})$ to an affine twist $(\tau\otimes\tau)({\cal F})$.
Let us denote $(\tau\otimes\tau)({\cal F}_{CG}')$ by $\hat{{\cal
F}}_{CG}^{(\tau)}$.
\begin{theorem}
Let $\hat{{\cal R}}^{'}$ be the universal $R$-matrix for
$U_{q}(\hat{\mathfrak{sl}}_{N})$. Then the $R$-matrix \\$\hat{{\cal
F}}_{CG}^{(\tau)21}\hat{{\cal R}}^{'}(\hat{{\cal F}}_{CG}^{(\tau)})^{-1}$
quantizes the quasi-trigonometric $r$-matrix obtained in Proposition \ref{p3}.
\end{theorem}
Both theorems can be proved straightforward.

\subsection{$\omega$-affinization and quantization of rational $r$-matrices}

Aim of this section is to quantize certain rational $r$-matrices. Of course, we
would like to use rational degeneration of the affine twists constructed above.
Unfortunately, we cannot do this directly because in both cases ${\rm
lim}_{q\rightarrow 1}{\cal F}=1\otimes 1$. Therefore, propose a new method
which we call $\omega$-affinization. We begin with the following result.
\begin{theorem}
Let $\pi:U_{q}(\mathfrak{sl}_{N}[u])\rightarrow U_{q}(\mathfrak{sl}_{N})$ be
the canonical projection sending all the affine generators to zero. Let
$\Psi\in U_{q}(\mathfrak{sl}_{N}[u])\otimes U_{q}(\hat{\mathfrak{sl}}_{N})$ be
invertible and such that
\begin{equation}\label{rat1}
\Psi=\Psi_{1}\Psi_{2}~,
\end{equation}
where $\Psi_{2}=(\pi\otimes{\rm id})(\Psi)$. Further, let $\Psi_{2}$ be a twist
and let the following two relations hold:
\begin{equation}\label{rat2}
\Psi_{1}^{23} \Psi_{2}^{12}\;=\;\Psi_{2}^{12}\Psi_{1}^{23},\qquad (
\pi\otimes{\rm id})(\Delta\otimes{\rm id})(\Psi_{1})\;=\; \Psi_{1}^{23}
\end{equation}
Finally, let there exist $\omega\in U_{q}(\mathfrak{sl}_{N}[u])[[\zeta]]$ such
that
\begin{equation}\label{rat3}
\Psi_{\omega}\,:=\;(\omega\otimes \omega) \Psi \Delta(\omega^{-1})\in
U_{q}(\mathfrak{sl}_{N})\otimes U_{q}(\hat{\mathfrak{sl}}_{N})~.
\end{equation}
Then $\Psi$ is a twist.
\label{affinizator}
\end{theorem}
\begin{proof}
Let us consider the Drinfeld associator
\begin{equation}\label{rat4}
{\rm Assoc}(\Psi)=\Psi^{12}(\Delta\otimes{\rm id})(\Psi) ({\rm
id}\otimes\Delta)(\Psi^{-1})(\Psi^{-1})^{23}
\end{equation}
and the following one equivalent to it
\begin{equation}\label{rat5}
\begin{array}{rcl}
{\rm Assoc}(\Psi_{\omega})\!\!&=\!\!&(\omega\otimes\omega \otimes\omega)({\rm
Assoc}(\Psi))(\omega^{-1}\otimes\omega^{-1}\otimes\omega^{-1})
\\[7pt]
\!\!&=\!\!& \Psi^{12}_{\omega}(\Delta\otimes{\rm id})(\Psi_{\omega})({\rm
id}\otimes\Delta)(\Psi^{-1}_{\omega})(\Psi^{-1}_{\omega})^{23}.
\end{array}
\end{equation}\label{rat6}
By (\ref{rat3}) we have
\begin{equation}
{\rm Assoc}(\Psi_{\omega})=(\pi\otimes{\rm id}\otimes{\rm id}) ({\rm
Assoc}(\Psi_{\omega}))~.
\end{equation}
Further, we take into account the following considerations:
\begin{equation}\label{rat7}
(\pi\otimes{\rm id})(\Psi_{2})\;=\;(\pi\otimes{\rm id})(\pi\otimes{\rm id})
\bigl(\Psi_{1}\Psi_{2}\bigr)\;=\;(\pi\otimes{\rm id})(\Psi_{1})(\pi\otimes{\rm
id})(\Psi_{2})~,
\end{equation}
what implies that $(\pi\otimes{\rm id})(\Psi_{1})=1\otimes1$. Moreover, the
latter also implies that
\begin{equation}\label{rat8}
(\pi\otimes{\rm id}\otimes{\rm id})({\rm id}\otimes
\Delta)(\Psi_{1}^{-1})\;=\;({\rm id}\otimes \Delta)(\pi\otimes{\rm
id})(\Psi_{1}^{-1})\;=\;1\otimes1\otimes1~.
\end{equation}
Thus,
\begin{equation}\label{rat9}
\begin{array}{rcl}
&&{\rm Assoc}(\Psi_{\omega})=
\\[5pt]
&&\quad={\rm
Ad}(\pi(\omega)\otimes\omega\otimes\omega)\bigr(\Psi_{2}^{12}\Psi^{23}_{1}
(\Delta\otimes{\rm id})(\Psi_{2})({\rm id}\otimes\Delta)(\Psi_{2}^{-1})
(\Psi_{2}^{-1})^{23}(\Psi_{1}^{-1})^{23}\bigr)
\\[7pt]
&&\quad={\rm Ad}(\pi(\omega)\otimes \omega\otimes\omega)\bigl(\Psi^{23}_{1}{\rm
Assoc}(\Psi_{2})(\Psi^{23}_{1})^{-1}\bigr)~.
\end{array}
\end{equation}
Since $\Psi_{2}$ is a twist  we deduce that ${\rm Assoc}(\Psi_{\omega})=
1\otimes 1\otimes 1$ and therefore ${\rm Assoc}(\Psi)=1\otimes 1\otimes 1$ what
proves the theorem.
\end{proof}
An element $\omega$ satifying conditions of Theorem \ref{affinizator} we will
call affinizator.
\begin{cor}
Let $\Psi_{1},\mathop{}\Psi_{2},\mathop{}\omega$ satisfy the conditions of
Theorem \ref{affinizator}. Then
\begin{equation}\label{rat10}
\Psi_{1}=(\omega^{-1}\pi(\omega)\otimes{\rm id}) \Psi_{2}\left\{(\pi\otimes{\rm
id})\Delta(\omega^{-1})\right\} \Delta(\omega)\Psi_{2}^{-1}.
\end{equation}
\label{psi}
\end{cor}
\begin{proof}
We have
\begin{equation}\label{rat11}
\Psi_{\omega}\;=\;(\pi\otimes{\rm id})(\Psi_{\omega})\;=\;
(\pi(\omega)\otimes\omega)\Psi_{2}(\pi\otimes{\rm id}) \Delta(\omega^{-1})~,
\end{equation}
because $\Psi_{2}=(\pi\otimes{\rm id})(\Psi)$. Now we see that
\begin{equation}\label{rat12}
(\omega\otimes\omega)(\Psi_{1}\Psi_{2})\Delta(\omega^{-1})\;=\;
(\pi(\omega)\otimes\omega)\Psi_{2}(\pi\otimes{\rm id})\Delta(\omega^{-1})~,
\end{equation}
what yields the required expression for $\Psi_{1}$.
\end{proof}

Now we would like to explain how $\omega$-affinization can be used to find a
Yangian degeneration of the affine Cremmer--Gervais twists. Let us consider the
case $\mathfrak{sl}_3$. We set
\begin{equation}\label{rat13}
\Psi_{2}:={\cal F}_{CG_{3}}^{(\tau)}:=\exp_{q^{2}}(-(q-q^{-1})
\zeta\mathop{}\hat{e}_{12}^{(0)}\otimes \hat{e}_{32}^{(0)})\cdot\hat{{\cal
K}_3}~,
\end{equation}
where
\begin{equation}\label{rat14}
\hat{\cal K}_3=q^{\frac 49 h_{12}\otimes h_{12}+\frac 29 h_{12}\otimes
h_{23}+\frac 59 h_{23}\otimes h_{12}+ \frac 79 h_{23}\otimes h_{23}}
\end{equation}
with $h_{ij}:=e_{ii}-e_{jj}$. The twisting two-tensor (\ref{rat11}) belongs to
$U_{q}(\mathfrak{sl}_{3})\otimes U_{q}(\mathfrak{sl}_{3})[[\zeta]]$.

The following affinizator $\omega_{3}^{\rm long}$ was constructed in
\cite{Sam1}. It is given by the following formula
\begin{equation}\label{rat15}
\begin{array}{rcl}
\omega_{3}^{\rm long}\!\!&=\!\!&
\displaystyle\exp_{q^{2}}(\frac{\zeta}{1-q^{2}}\mathop{}q^{2h_{\alpha}^{\perp}}
\hat{e}_{21}^{(1)})\exp_{q^{2}}(-\frac{q\zeta^{2}}{(1-q^{2})^{2}}
q^{2h_{\beta}^{\perp}}\hat{e}_{31}^{(1)})
\\[9pt]
&&\times\,\displaystyle\exp_{q^{-2}}(\frac{\zeta^{2}}{1-q^{2}}
\mathop{}\hat{e}_{32}^{(0)})\exp_{q^{-2}}(\frac{\zeta}{1-q^{2}}\mathop{}
\hat{e}_{21}^{(0)})\exp_{q^{-2}}(\frac{\zeta^{2}}{1-q^{2}}\mathop{}
\hat{e}_{32}^{(0)})~,
\end{array}
\end{equation}
where $h_{\alpha}^{\perp}=\frac 13 (e_{11}+e_{22})-\frac 23 e_{33}$ and
$h_{\alpha}^{\perp}=\frac 23 e_{11}- \frac 13 (e_{22}+e_{33})$.

For convenience sake we remind the reader that
\begin{equation}\label{rat16}
\begin{array}{lcl}
\hat{e}_{12}^{(0)}=e_{12}^{0}q^{h_{\beta}^{\perp}-h_{\alpha}^{\perp}}&&
\hat{e}_{21}^{(0)}=q^{h_{\beta}^{\perp}}e_{21}^{(0)}\\[2ex]
\hat{e}_{32}^{(0)}=q^{-h_{\beta}^{\perp}}e_{32}^{(0)},&&
\hat{e}_{31}^{0}=e_{32}^{0}e_{21}^{(0)}-q^{-1}e_{21}^{(0)}e_{32}^{(0)}\\[2ex]
\hat{e}_{31}^{(1)}=q^{h_{\alpha}^{\perp}-h_{\beta}^{\perp}}e_{31}^{(1)},&&
\hat{e}_{32}^{(1)}=e_{12}^{(0)}e_{31}^{(1)}-q e_{31}^{(1)}e_{12}^{(0)}\\[2ex]
\end{array}
\end{equation}
\begin{theorem}
The elements $\omega_{3}^{\rm long}$, $\Psi_{2}={\cal F}_{CG_{3}}^{(\tau)}$ and
$\Psi_{1}$ constructed according to Corollary \ref{psi} satisfy the conditions
of Theorem \ref{affinizator} and consequently $\Psi_{\omega}=(\pi\otimes{\rm
id})(\omega\otimes\omega)\Psi_{2}\Delta(\omega^{-1})$ is a twist.
\end{theorem}
\begin{proof}
Straigtforward.
\end{proof}

It turns out that $\Psi_{\omega}$ has a rational degeneration. To define this
rational degeneration we have to introduce the so-called $f$-generators:
\begin{equation}\label{rat17}
\begin{array}{lcl}
f_{0}=(q-q^{-1})\mathop{}\hat{e}_{31}^{(0)},&&
f_{1}=q^{2h_{\beta}^{\perp}}
\hat{e}_{31}^{(1)}+q^{-1}\zeta\mathop{}\hat{e}_{31}^{(0)},\\[2ex]
f_{2}=(1-q^{-2})\mathop{}\hat{e}_{32}^{(0)},&& f_{3}=q^{h_{\alpha}^{\perp}}
\hat{e}_{32}^{(1)}-\zeta\mathop{}\hat{e}_{32}^{(0)}.
\end{array}
\end{equation}
Let us consider the Hopf subalgebra of  $U_{q}^{\hat{{\cal
K}}_3}(\hat{\mathfrak{sl}}_{3})$ generated by
$$\{h_{12}, h_{23}, f_{0}, f_{1},f_{2}, f_{3},
\hat{e}_{12}^{(0)},\hat{e}_{21}^{(0)}\}.$$
When $q\rightarrow 1$ we obtain the following Yangian twist (see \cite{Sam1}):
\begin{equation}\label{rat18}
\begin{array}{rcl}
\overline{{\cal
F}}_{1}\!\!&=\!\!&(1\otimes1-\zeta\mathop{}1\otimes\overline{f}_{3}-
\zeta^{2}\mathop{}h_{\beta}^{\perp}\otimes
\overline{f}_{2})^{(-h_{\beta}^{\perp}\otimes 1)}(1\otimes 1+
\zeta\mathop{}1\otimes\overline{e^{(0)}_{21}})
^{(-h_{\beta}^{\perp}\otimes 1)}
\\[7pt]
&&\times\;\exp(\zeta^{2}\mathop{}\overline{e^{(0)}_{12}}h_{13}
\otimes\overline{f}_{0})\exp(-\zeta\mathop{}
\overline{e^{(0)}_{12}}\otimes\overline{f}_{1})\cdot
\exp(-\zeta\mathop{}\overline{{e}_{12}^{(0)}}\otimes
\overline{f}_{2})
\\[7pt]
&&\times\; (1\otimes
1-\zeta\mathop{}1\otimes\overline{f}_{3}-\zeta^{2}\mathop{}
h_{\alpha}^{\perp}\otimes\overline{f}_{2})^{((h_{\beta}^{\perp}-
h_{\alpha}^{\perp})\otimes1)}~,
\end{array}
\end{equation}
where the overlined generators are the generators of $Y(\mathfrak{sl}_{3})$. In
the evaluation representation we have:
\begin{equation}\label{rat19}
\begin{array}{ccc}
\overline{f}_{0}\mapsto e_{31}~,&& \overline{f}_{1}\mapsto u\mathop{}e_{31}
\\[5pt]
\overline{f}_{2}\mapsto e_{32}~, &&
{f}_{3}\mapsto u\mathop{}e_{32}
\\[5pt]
\overline{\vphantom{f}e}_{21}\mapsto e_{21}~,&&
\overline{\vphantom{f}e}_{12}\mapsto e_{12}~.
\end{array}
\end{equation}
Therefore we have obtained the following result:
\begin{theorem}
The Yangian twist $\overline{{\cal F}}_{1}$ quantizes the following classical
rational $r-$matrix
\begin{equation}\label{rat20}
\begin{array}{rcl}
r(u,v)\!\!&=\!\!&\displaystyle\frac{\Omega}{u-v}+h_{\alpha}^{\perp}\otimes
ve_{32}-ue_{32} \otimes h_{\alpha}^{\perp}+h_{\beta}^{\perp}\wedge e_{21}
\\[10pt]
\!\!&&+e_{21}\otimes ve_{31}-ue_{31}\otimes e_{21}+e_{12}\wedge e_{32}~.
\end{array}
\end{equation}
\end{theorem}

To obtain a quantization of the second non-trivial rational $r$-matrix for
$\mathfrak{sl}_{3}$ we take the following affinizator
$\omega_{3}^{\rm\mathop{}short}$ and apply it to $\Psi_{2}=q^{r_0(3)}$, where
the Cartan part of the Cremmer-Gervais constant $r$-matrix for
$\mathfrak{sl}_{3}$ has the form:
\begin{equation}\label{rat21}
r_0(3)\;=\;\frac{2}{3}\bigr(h_{\alpha_1}\otimes h_{\alpha_1}+
h_{\alpha_2}\otimes h_{\alpha_2}\bigl)+\frac{1}{3}\bigr(h_{\alpha_1}\otimes
h_{\alpha_2}+h_{\alpha_2}\otimes
h_{\alpha_1}\bigl)+\frac{1}{6}h_{\alpha_1}\wedge h_{\alpha_2}~.
\end{equation}
We have
\begin{equation}\label{rat22}
\omega_{3}^{\rm short}\;=\;\exp_{q^{-2}}\bigl(\zeta\mathop{}\hat{e}_{21}^{(0)}
\bigr)\,\exp_{q^{2}}\Bigl(-\frac{\zeta}{1-q^{2}}\mathop{}q^{2h_{\alpha}^{\perp}}
\mathop{}\hat{e}_{31}^{(1)}\Bigr)\,\exp_{q^{-2}}\Bigl(\frac{\zeta}{1-q^{2}}\mathop{}
\hat{e}_{32}^{(0)}\Bigr)~,
\end{equation}
where
$$
\begin{array}{ccccc}
\hat{e}_{21}^{(0)}=q^{-\frac 13(h_{12}-h_{23})} e_{21}^{(0)},&&
\hat{e}_{32}^{(0)}=q^{-h_{\alpha}^{\perp}}e_{32}^{(0)},&&
\hat{e}_{31}^{(1)}=q^{-h_{\alpha}^{\perp}}e_{31}^{(1)}~.
\end{array}
$$
We have to calculate
\begin{equation}\label{rat23}
{\rm Aff}_{\omega^{\rm\mathop{} short}_3}(q^{r_0(3)})\;:= \;(\pi\otimes{\rm
id})\circ\Bigl((\omega_{3}^{\rm short}\otimes \omega_{3}^{\rm short})
q^{r_0(3)}\Delta(\omega_{3}^{\rm\ short})^{-1}\Bigr)~.
\end{equation}

Using standard commutation relations between $q$-exponents, the formula
(\ref{rat4}) can be brought to the following form:
\begin{equation}\label{rat24}
\begin{array}{rcl}
&&\Bigr(1\otimes 1+\zeta\mathop{}1\otimes q^{2h_{\alpha}^{\perp}}
\hat{e}_{31}^{(1)}+\zeta\mathop{}q^{-2h_{\alpha}^{\perp}}\otimes\bigr({\rm
Ad}\exp_{q^{2}}(\zeta\hat{e}_{21}^{(0)})\bigr)(\hat{e}_{32}^{(0)})\Bigr)_{q^{2}}
^{(-h_{\alpha_{1}}^{\perp}\otimes 1)}
\\[7pt]
&&\qquad\times\;\Bigl(1\otimes 1+\zeta (1-q^{2})\mathop{}
1\otimes\hat{e}_{21}^{(0)}\Bigr)_{q^{-2}}^{(-\frac 13 (h_{12}-h_{23})\otimes
1)}q^{r_0(3)}~.
\end{array}
\end{equation}
The $q$-Hadamard formula allows us to calculate the Ad-term explicitly:
\begin{equation}\label{rat25}
\bigl({\rm Ad}\exp_{q^{-2}}(\zeta\mathop{}\hat{e}_{21}^{(0)})\bigr)
(\hat{e}_{21}^{(0)})= \hat{e}_{21}^{(0)}+\zeta\mathop{}q^{-h_{\beta}^{\perp}}
e_{31}^{(0)}~,
\end{equation}
where $e_{31}^{(0)}:= e_{21}^{(0)}e_{32}^{(0)}-q\mathop{}
e_{32}^{(0)}e_{21}^{(0)}$. To define a rational degeneration we introduce
$g$-generators, which satisfy the Yangian relations as $q\to 1$:
$$
\begin{array}{ccccc}
g_{0}=(q-q^{-1})q^{-h_{\beta}^{\perp}}e^{(0)}_{31},&&
g_{1}=q^{2h_{\alpha}^{\perp}}\hat{e}_{31}^{(1)}+
\zeta\mathop{}q^{-h_{\beta}^{\perp}}
e^{(0)}_{31},&&
g_{2}=(q^{2}-1)\mathop{}\hat{e}_{21}^{(0)}.
\end{array}
$$
Using $g$-generators we can calculate the rational degeneration of the twist
${\rm Aff}_{\omega^{\rm\mathop{}short}_3}(q^{r_0(3)})$:
\begin{equation}\label{rat26}
\begin{array}{rcl}
\overline{{\cal F}}_{2}\!\!&=\!\!&\Bigl(1\otimes
1+\zeta\mathop{}1\otimes(\overline{g}_{1}+\overline{e^{(0)}_{32}})
-\zeta^2\mathop{}h_{\alpha}^{\perp}\otimes\overline{g}_{0}\Bigr)
^{(-h_{\alpha}^{\perp}\otimes 1)}
\\[7pt]&&\times\;\bigl(1\otimes 1-\zeta\mathop{}1\otimes
\overline{g}_{2}\bigr)^{(-\frac 13(h_{12}-h_{23})\otimes 1)}.
\end{array}
\end{equation}
\begin{theorem}
This Yangian twist $\overline{{\cal F}}_{2}$ quantizes the following rational
$r$-matrix:
\begin{equation}\label{rat27}
r(u,v)\;=\;\frac{\Omega}{u-v}-u\mathop{}e_{31}\otimes
h_{\alpha}^{\perp}+v\mathop{}h_{\alpha}^{\perp}\otimes e_{31}
+h_{\alpha}^{\perp}\wedge e_{32}-\frac 13 (h_{12}-h_{23})\wedge e_{21}~.
\end{equation}
\end{theorem}
Therefore we have quantized all non-trivial rational $r$-matrices for
$\mathfrak{sl}_3$ classified in \cite{S}.

\subsection*{Acknowledgments}
The paper has been partially supported by the Royal Swedish Academy of Sciences
under the program "Cooperation between Sweden and former USSR" and the grants
RFBR-05-01-01086, INTAS-OPEN-03-51-3350 (V.N.T.).


\begin{thebibliography}{2}

\bibitem{CM} A. Connes and H. Moscovici,
\textsl{Rankin--Cohen Brackets and the Hopf Algebra of Transverse Geometry},
Moscow Math. J., \textbf{4}(1) (2004), 111--130, 311.

\bibitem{BD} A. Belavin and V. Drinfeld,
\textsl{Solutions of the classical Yang--Baxter equation for simple Lie algebras},
Functional Anal. Appl., \textbf{16}(3) (1983), 159--180; {\it translated from}
Funktsional. Anal. i Prilozhen, \textbf{16} (1982), 1--29 (Russian).

\bibitem{CG} E. Cremmer and J.L. Gervais,
\textsl{The quantum group structure associated with nonlineary extended Virasoro
algebras}, Comm. Math. Phys., \textbf{134}(3) (1990), 619--632.

\bibitem{DK} V. Dergachev and A.A. Kirillov,
\textsl{Index of Lie algebras of seaweed type}, J. Lie Theory,
\textbf{10}(2) (2000), 331--343.

\bibitem{GG} M. Gerstenhaber, A. Giaquinto,
\textsl{Boundary solutions of the classical Yang-Baxter equation},
Lett. Math. Phys. \textbf{40}(4) (1997), 337--353.

\bibitem{D} V.G. Drinfeld,
\textsl{On some unsolved problems in quantum group theory}, in: Quantum
groups, (Springer, Berlin, 1992), 
Lecture Notes in Math., \textbf{1510}, 1--8.

\bibitem{ESS} P. Etingof, T. Schedler, and Schiffmann,
\textsl{Explicit quantization of dynamical $r$-matrices for finite dimensional
semisimple Lie algebras}, J. Amer. Math. Soc., \textbf{13} (2000), 595--609.

\bibitem{H} T.J. Hodges,
\textsl{Nonstandard quantum groups associated to certain Belavin--Drinfeld
triples}. Perspectives on quantization (South Hadley, MA, 1996), 63--70,
Contemp. Math., \textbf{214}, Amer. Math. Soc., Providence, RI, 1998.

\bibitem{IO} A.P. Isaev and O.V. Ogievetsky,R
\textsl{On Quantization of $r$-matrices for Belavin-Drinfeld Triples},
Phys. Atomic Nuclei, \textbf{64}(12) (2001), 2126--2130, {\tt math.QA/0010190}.

\bibitem{KPST} S.M. Khoroshkin, I.I. Pop, A.A. Stolin, and Tolstoy,
\textsl{On some Lie bialgebra structures on polynomial algebras and their
quantization}, preprint, Mittag-Leffler Institute, Sweden, 2004.

\bibitem{KST} S.M. Khoroshkin, A.A. Stolin, and Tolstoy,
\textsl{Deformation of Yangian $Y({\rm sl}\sb 2)$},
Comm. Algebra, \textbf{26}(4)  (1998), 1041--1055.

\bibitem{KST1} S.M. Khoroshkin, A.A. Stolin, and V.N. Tolstoy,
\textsl{$q$-Power function over $q$-commuting variables and deformed $XXX$ and
$XXZ$ chains}, Phys. Atomic Nuclei, \textbf{64}(12) (2001), 2173--2178;
{\it translated from} Yad. Fiz. \textbf{64}(12) (2001), 2262--2267.

\bibitem{KT1} S.M. Khoroshkin and V. Tolstoy,
\textsl{Universal $R$-matrix for quantized {super}algebras}, Comm. Math. Phys.
\textbf{141}(3) (1991), 599--617.

\bibitem{KT2} S.M. Khoroshkin and V.N. Tolstoy,
\textsl{Twisting of quantum (super)algebras. Connection of Drinfeld's and
Cartan--Weyl realizations for quantum affine algebras}, MPIM preprint, MPI/94-23,
pp. 1--29 (Bonn, 1994); {\tt arXiv:hep-th/9404036}.

\bibitem{KM} P.P Kulish and A.I. Mudrov,
\textsl{Universal $R$-matrix for esoteric quantum groups},
Lett. Math. Phys. \textbf{47} (2) (1999), 139--148.

\bibitem{P} D.I.Panyushev,
\textsl{Inductive formulas for the index of seaweed Lie algebras},
Mosc. Math. J., \textbf{1}(2) (2001), 221--241.

\bibitem{Sam1} M. Samsonov,
\textsl{Semi-classical Twists for $\mathfrak{sl}_{3}$ and $\mathfrak{sl}_{4}$
Boundary $r$-matrices of Cremmer-Gervais type}, Lett. Math. Phys.,
\textbf{72}(3) (2005), 197--210.

\bibitem{Sam2} M. Samsonov,
\textsl{Quantization of semi-classical twists and noncommutative geometry},
Lett. Math. Phys., \textbf{75}(1) (2006), 63--77, {\tt math.QA/0309311}.

\bibitem{S} A. Stolin,
\textsl{On rational solutions of Yang--Baxter equation for ${\mathfrak{sl}}(n)$},
Math. Scand., \textbf{69}(1) (1991), 57--80.

\bibitem{St} A. Stolin,
\textsl{Some remarks on Lie bialgebra structures on simple complex Lie algebras},
Comm. Algebra, \textbf{27}(9) (1999), 4289--4302.

\bibitem{T1} V.N. Tolstoy,
\textsl{Extremal projectors for quantized Kac--Moody superalgebras
and some of their applications}. Lecture Notes in Phys., (Springer, Berlin),
\textbf{370} (1990), 118--125.

\bibitem{TK} V.N. Tolstoy and S.M. Khoroshkin,
\textsl{Universal $R$-matrix for quantized nontwisted affine Lie algebras},
Func. Anal. Appl. \textbf{26}(1) (1992), 69--71; {\it translated from}
Funktsional. Anal. i Prilozhen, \textbf{26}(1) (1992), 85--88 (Russian).

\bibitem{T2} V.N. Tolstoy,
\textsl{From quantum affine Kac--Moody algebra to Drinfeldians and
Yangians}, in: Kac-Moody Lie algebras and related topics,
AMS 2004, Contemporary Mathematics, CONM/\textbf{343},
349--370; {\tt arXiv:math.QA/0212370}.

\bibitem{T3} V.N. Tolstoy,
\textsl{Super-Drinfeldian and super-Yangian for the superalgebra $U_q(sl(n|m))$},
Phys. Atomic Nuclei, \textbf{64}(12) (2001), 2179--2184;
{\it translated from} Yad. Fiz. \textbf{64}(12) (2001), 2268--2273.

\end{thebibliography}
\end{document}